\documentclass{amsart}

\usepackage{color}
\usepackage{mathscinet}
\usepackage{epigraph}
\newtheorem{theorem}{Theorem}[section]

\theoremstyle{definition}

\theoremstyle{remark}

\numberwithin{equation}{section}

\begin{document}

\title{Geometry and physics: an overview}


\author{David R. Morrison}
\address{Departments of Mathematics and Physics, University of California, 
  Santa Barbara CA 93106}
\email{drm@math.ucsb.edu}
\thanks{I am grateful to Arun Debray and Sean Pohorence for sharing the notes
they took during the lecture, and to Andy Neitzke for a very useful remark.  
I am also grateful to Greg Moore and Samson Shatashvili for their comments on an earlier draft of this article.
This research was 
partially supported by 
NSF grant PHY-1620842.} 




\date{}

\begin{abstract}
We present some episodes from the history of interactions
between geometry and physics over the past century.
\end{abstract}

\maketitle



I was asked by the conference organizers to survey the modern interactions
between geometry, topology, and physics in a one-hour lecture.  Rather than
attempting to be comprehensive, I  have chosen
 to provide some vignettes drawn from recent history which
emphasize the impact 
each field has had upon the other.  As part of the story, 
I will trace a particular circle of ideas from physics to math, 
back to physics,
and on to math once again, and show the significant impacts at each stage along
the way.  My general theme is gauge theory, geometry, and topology.

\section{Dirac quantization} \label{sec:dirac}

\epigraph{The most powerful method of advance that can be suggested at 
present is to employ all
the resources of pure mathematics in attempts to perfect and generalize 
the mathematical
formalism that forms the existing basis of theoretical physics, and after 
each success in this
direction, to try to interpret the new mathematical features in terms of 
physical entities.}{Paul Dirac (1931)}

Paul Dirac laid out the
manifesto quoted above  for the proper
interaction between mathematics
and physics in the course of
a beautiful 1931 paper \cite{dirac-quantised} devoted to
the theoretical investigation of magnetic monopoles.

It had long been observed that Maxwell's equations for electromagnetism
can be made symmetric between electricity and magnetism by introducing
particles carrying a net magnetic charge: the so-called 
magnetic monopoles.\footnote{I
hasten to point out that magnetic monopoles have not (yet)
been observed in nature.}  
The quantum version of Maxwell's theory, however, relies on the 
existence of an electromagnetic potential and prior to Dirac's work it
was believed that monopoles would prevent such a potential from being
defined.  Dirac showed that by using a careful interpretation of the role
of the potential, quantum physics could be formulated using it (in spite
of the definitional problem).  Morever,
there was a surprise: the existence of a magnetic monopole forces the ``charge
quantization'' of electrically charged particles such as the electron.

There was some interesting topology in Dirac's original argument,
a variant of which I will now review.
Consider  the wave function of an electron 
in the presence of a magnetic
monopole of magnetic charge $g$, 
located at the origin in $\mathbb{R}^3$, and consider a
vector potential $A$ for the magnetic field.  One might suppose that the
vector potential could be defined throughout $\mathbb{R}^3-\{\vec{0}\}$,
but that is not possible for the following reason.  If we consider a
sphere of radius $r$ and a circle $\gamma_\theta$ at fixed 
spherical polar coordinate
angle $\theta$ which bounds a spherical cap $\Sigma_\theta$ containing the
north pole, then the line integral
\begin{equation}
\int_{\gamma_\theta} A \, d\varphi 
= \int_{\Sigma_\theta} \nabla \times A \, d\sigma
= \int_{\Sigma_\theta} B \, d\sigma
\end{equation}
calculates the magnetic flux through $\Sigma_\theta$ (by Stokes' theorem).
That flux, which can be calculated to be $2\pi g (1-\cos\theta)$ increases
continuously as $\theta$ ranges from $0$ to $\pi$.  But at $\theta=\pi$
the circle has shrunk to zero size so that the integral must vanish,
which is a contradiction!

Dirac's interpretation was that the vector potential must become singular
along some semi-infinite string anchored at the monopole.  If we direct
such a string from the origin through the south pole of the sphere,
the contradiction is removed.  The integral around a zero-size circle 
at the south pole
cannot be computed directly due to the singularity in the vector potential,
but takes the value $4\pi g$ thanks to the integral over the sphere.

The magnetic field causes a change in the phase of the wave function
of an electric particle around a circle.
For an
electric particle of charge $e$,
if the change in phase of the wave function around a circle 
is zero at the north pole, and it must take the value
\begin{equation}
\frac e{\hbar c} \int_{S^2} B\, d\sigma= \frac{4\pi ge}{\hbar c}
\end{equation}
at the south pole.
In a quantum theory, 
a change of phase
in the wave function is physically indetectable if it is a multiple of $2\pi$.
The conclusion is that in order for the quantum theory to be well-defined,
we must have
\begin{equation} \frac{4\pi ge}{\hbar c} = 2\pi N \end{equation}
for some integer $N$.  In other words, if there is a single monopole of
magnetic charge $g$, then all electric charges must be integer multiples
of $\hbar c/2g$.  This is Dirac's famous quantization condition.

Conversely, since we know experimentally the smallest possible electric
charge $e_0$, all magnetic charges must be integer multiples of
$\hbar c/2e_0$.

In modern topological terms,\footnote{This mathematical
argument first appeared  in
the physics literature \cite{MR0426712}.}  the vector potential $A$ should
only be locally defined, with a change between northern hemisphere and
southern hemisphere specified by a change of potential in a neighborhood
of the equator, which amounts to a change by a variable element of $U(1)$
(known as an {\em abelian gauge transformation}).The vector potential can be regarded as a globally defined object if
it is treated as a section of a bundle over $\mathbb{R}^3-\{\vec{0}\}$
(or over $S^2$) whose transition functions are given by that variable
element.

The abelian gauge transformation is specified by a map
$S^1\to U(1)$. Since $\pi_1(U(1))\cong \mathbb{Z}$, the topological
type of the bundle is determined by an integer (the first
Chern class of the bundle).  Physically, this integer is identified
with the multiple of $\hbar c/2e_0$ which gives the magnetic charge $g$
of the monopole in question.

Dirac did not formulate his result in terms of a bundle, but his
treatment is surprisingly modern, given the relative novelty
of topological concepts at that time.  This result turns out to be
only the first step in an important series of interactions between
topology, geometry, and physics.

\section{Missed opportunities} \label{sec:divorce}

\epigraph{As a working physicist, I am acutely aware of the fact that the
marriage between mathematics and physics, which was so fruitful in past
centuries, has recently ended in divorce.}{Freeman Dyson (1972)}

In his 1972 Gibbs lecture \cite{MR0522147}, Freeman Dyson lamented
the divide between mathematics and physics which existed at that time,
drawing a number of examples from history.  One of his principal
examples centered around Feynman's approach to
 the study of relativistic quantum field theory.

The years immediately following World War II saw rapid advances in
the study of relativistic quantum field theory, led by Dyson, Feynman,
Schwinger, and Tomonaga.  Feynman's approach was based on his
``sum over histories'' idea.  In the Lagrangian formulation of classical
physics, the evolution of a physical system from a starting time to
an ending time is the one which minimizes the physical {\em action}\/
of the system.  Such a minimum can be determined using the techniques
of the calculus of variations to find the path through the configuration
space which extremizes the action.  In Feynman's
``sum over histories'' approach to quantum physics, 
{\em all}\/ paths must be considered,
and the probability that a particular path is followed is proportional to
the exponential of the negative of the action.  Thus, the probability
will be highest along the classical path, but the quantum theory
requires consideration of contributions from other paths.

To determine physically measurable quantities, it is necessary to
integrate over the space of all paths, and this is where the mathematical 
trouble
arises: it is unknown how to carry out such integrals.  One of the
difficulties is a choice of measure on the space of paths which determines
how different paths are to be weighted.

Nevertheless, by considering theories which are perturbations of ``free''
theories, Feynman was able to give a prescription for an asymptotic
series describing the (purported) answer to the path integral, as an
infinite sum indexed by
the famous {\em Feynman diagrams.}  The mathematical consistency of
the path integral values as determined by Feynman diagrams has never
been established, but in the hands of skilled practitioners unique
answers are produced.  Remarkably, when this approach is used to study
quantum electrodynamics, it is incredibly precise:  for example,
the best theoretical and experimental values of the anomalous magnetic
dipole moment of the electron agree to at least eight significant figures
(see \cite[Chapter 6]{MR1402248}).

As Dyson remarks, Feynman himself was not concerned with mathematical
rigor.  The fact that physicists were obtaining such spectacular agreement
between theory and experiment likely contributed to a similar attitude
among many physicists from that era:  they did not need modern mathematics
to accomplish their goals.  And so Dirac's manifesto languished.

\section{Yang--Mills theory and connections on fiber bundles} \label{sec:YM}

\epigraph{\dots we are concerned with the necessary {\em concepts to
describe the physics of gauge theories.}  It is remarkable that
these concepts have already been studied as mathematical 
constructs.}{Tai Tsun Wu and Chen Ning Yang (1975)}

In 1954, during the era of minimal communication between mathematics
and theoretical physics, C.~N. Yang and R.~L. Mills 
\cite{MR0065437} introduced  {\em gauge
transformations}\/ consisting of locally varying symmetries taking
values in a compact Lie group\footnote{To be precise,
\cite{MR0065437} treats the case $G=SU(2)$ but the notion was
soon generalized to an arbitrary compact group.} $G$, and studied
physical theories which are invariant under such gauge transformations.
These generalized the already-familiar
{\em abelian gauge transformations}\/ from electromagnetism -- the same
ones we encountered in Section~\ref{sec:dirac} -- for which $G=U(1)$.
These gauge theories (or ``Yang--Mills theories'') eventually became the 
basis of the Standard Model of particle physics, the formulation
of which was
finalized in the mid 1970s using the 
group\footnote{For a mathematical account of the standard model gauge
group and the connection between its representations and 
the elementary particles, see \cite{MR2651086}.}
$G=(SU(3)\times SU(2)\times U(1))/\mathbb{Z}_6$.

In the late 1960s and early 1970s, Yang
got acquainted with James
Simons, then the mathematics department chair at SUNY Stony Brook 
where Yang was a professor of physics.
In the course of their conversations,\footnote{Yang and Simons shared
memories of this period in a joint interview in 2008 \cite{simons-yang}.}
Yang and Simons came to recognize that
there were important similarities between formulas which were showing
up in Yang's work, and formulas which appeared in parts of 
mathematics which Simons was familiar with.  Simons identified the
relevant mathematics as the mathematical theory of {\em 
connections on fiber bundles},
and recommended that Yang consult Steenrod's foundational book on the subject 
\cite{MR0039258} (which coincidentally was published  just a few
years prior to the work of Yang and Mills).  Yang found the book difficult to
read, but through further discussions with Simons and other mathematicians
(including S.-S. Chern) he
came to appreciate the power of the mathematical tools which fiber bundle
theory offered.  By 1975, Yang had co-authored a paper with T.\ T.\ Wu
\cite{MR0426712} (quoted at the head of this section) which applied those
methods to problems in physics.  Within their paper,
Wu and Yang provided a dictionary
between the parallel concepts in physics and mathematics,
allowing the application of topological and geometric techniques
to the study of Yang--Mills theory.

Simons communicated these newly uncovered connections with physics to
Isadore Singer at MIT who in turn discussed them with Michael Atiyah
of Cambridge University.  Similar observations were
made independently by others.  For example, in Moscow Ludwig Faddeev was inspired to give a geometric interpretation of Yang--Mills theory \cite{MR3748689} by the book 
\cite{MR0088015}. 

A new chapter in the interaction between
mathematics and physics was about to open.

\section{Unreasonable effectiveness} \label{sec:instantons}

\epigraph{\dots mathematical concepts turn up in entirely
unexpected connections.  Moreover, they often permit an unexpectedly
close and accurate description of the phenomena in 
these connections.}{Eugene Wigner (1960)}

A key example of Wigner's ``unreasonable effectiveness of mathematics''
principle \cite{wigner-unreasonable} occurred in the mid 1970s, not
long after Dyson's Gibbs lecture.
Thanks to the
 opening of communication between mathematicians such as Atiyah and
Singer on the one hand and the gauge theory community in physics on the
other hand, 
when Polyakov proposed in \cite{Polyakov:1975rs} the importance of studying 
instantons
in Yang--Mills theory, mathematicians were ready to assist in finding
such instantons.

A {\em Yang--Mills instanton}\/ is a solution on $\mathbb{R}^4$
to  the Euclidean version%
\footnote{One of the important techniques for studying 
quantum field theory considers these
Euclidean theories as a kind
of  analytic
continuation from ordinary time $t$ to ``imaginary time'' $it$, 
via a procedure known as
{\em Wick rotation.}  It is beyond the scope of this lecture to explain
why this is relevant to physics.}  
of the Yang--Mills equations for a compact Lie group $G$,
which are the variational equations for the norm-squared 
\begin{equation}
\|F\|^2 = \int_{\mathbb{R}^4} \operatorname{tr} (F \wedge \star F)
\end{equation}
of the curvature $F$ of a connection $A$ on a principal $G$-bundle,
where the Hodge star operator is used to define the norm.
Such a solution is only interesting to physicists
if it has a suitably controlled behavior far from the origin, and the
initial assumption made by mathematicians in studying the problem
is that%
\footnote{It was later shown
by Uhlenbeck \cite{MR648355} that this is automatically true if
$\|F\|^2$ is finite.}
the solution extends to $S^4$.

After some progress had been made on the problem in the physics community,
it was given a purely mathematical formulation by Atiyah, Hitchin,
and Singer \cite{MR0458424}, who established the dimension of the
space of solutions.  Not long thereafter the problem was solved in general
by Atiyah, Drinfeld, Hitchin and Manin \cite{MR598562} via what came to be
known as the ``ADHM construction.''  These papers
 used techniques -- the Penrose
twistor transform, and the algebraic geometry of vector bundles -- which were
then unknown to physicists.

Let me briefly explain the topological setting of the Yang--Mills 
instanton problem (already described in one of the earliest papers
\cite{MR0434183}), which can be viewed in two ways.  From the perspective
of a solution on $\mathbb{R}^4$, all $G$-bundles are trivial but
there is an asymptotic behavior of an instanton which is determined by
the behavior of the connection on
the $S^3$ at infinity.  There must be a gauge transformation on $S^3$,
i.e., a map $S^3\to G$,
which trivializes the connection there.  Thus, the topology of an 
instanton is measured by $\pi_3(G)$, which is  isomorphic to
$\mathbb{Z}$ for any compact semisimple Lie group.  The resulting integer $k$
is called the {\em instanton number}.\footnote{Changing the orientation
of the four-manifold changes the sign of the instanton number, and
there is some ambiguity in the literature about how this is defined.
We have attempted to be self-consistent in this paper.}

From the perspective of a solution on $S^4$, the map $S^3\to G$
specifies the bundle by giving gluing data along the equator.
The topological classication of principal $G$-bundles on $S^4$
is via $\pi_3(G)$ and determines the  instanton number as before.  
To see how this is related to the curvature of a connection, it is
convenient to remember that the Hodge star operator $\star$ on a 
Riemannian four-manifold
squares to the identity on $2$-forms.  Thus, $F$ can be decomposed 
into its {\em self-dual}\/ and {\em anti-self-dual}\/ parts:
\begin{equation}
F = F_{\text{sd}} + F_{\text{asd}}
\end{equation}
where $\star F_{\text{sd}} = F_{\text{sd}}$ and
$\star F_{\text{asd}} = - F_{\text{asd}}$.
We then have 
\begin{equation}
\begin{aligned} \| F \|^2 &= \int_{S^4} \operatorname{tr} ( F \wedge \star F) 
= \| F_{\text{sd}} \|^2 + \| F_{\text{asd}} \|^2 \\
8\pi^2 k &= \int_{S^4} \operatorname{tr}( F \wedge F ) 
= \| F_{\text{sd}} \|^2 - \| F_{\text{asd}} \|^2,
\end{aligned}
\end{equation}
reflecting the topology of the solution.  We also see from this that the
minimal action solutions must be either self-dual or anti-self-dual
depending on the sign of $k$.

I was fortunate enough to attend the {\em Loeb lectures}\/ delivered at Harvard
University by Michael Atiyah in the spring of 1978
in which he explained the ADHM construction.%
\footnote{Atiyah
subsequently published  notes based on the Loeb lectures as well as similar
lectures delivered at two other places
 \cite{MR554924}.} 
The lectures were held 
in the physics department, and probably constituted the largest 
 meeting that had  been held up until that time at Harvard
between mathematicians and 
physicists.\footnote{Harvard was and is strong in the more traditional areas of
mathematical physics, but the ADHM construction
 was something new and different and 
attracted a big audience from outside that community.}
My memory is that the audience was roughly half and half:  there were
large numbers of
mathematicians as well as large numbers of physicists.  
To an algebraic geometry
graduate student such as myself, it was an amazing experience to see
my quite abstract corner of mathematics applied to the ``real world'' of
theoretical physics.  The lectures were frustrating in one sense:  Atiyah
made the pedagogical choice of treating algebraic geometry as a ``black box''
for the purpose of the lectures, so we didn't get to hear about the
details of the algebraic geometry!  (I believe that some of the physicists
were frustrated by this as well, since they missed the opportunity to
learn about the algebraic geometry.)
However, we did learn about the twistor transform which was new material
to many of us.

\section{Anomalies and index theory} \label{sec:anomalies}

\epigraph{From a more mathematical standpoint, the study of anomalies
has elicited very interesting applications of index theory in quantum
field theory.  [\dots\  A]  powerful form of the Atiyah--Singer index
theorem (the index theorem for families of elliptic operators) has
been used to provide a global understanding of the non-Abelian anomaly
as well as the gravitational anomalies.}{Luis Alvarez-Gaum\'e (1986)}

Let me now turn to another chapter in the math--physics dialogue of the
1970s and 1980s: the calculation of anomalies in quantum field theories
and  quantum theories of gravity.  Although not directly related to our main story it is another signification interaction of the era.

Classical physical theories are determined by the equations of
motion of the theory, but quantum theories require an understanding
of the physical ``action'' on a broad configuration space of
possible physical fields. Moreover,  Feynman's formulation 
requires an understanding of a measure on the space of paths through
that configuration space as well.  

In the 1960s, while attempting to understand
how symmetries of  classical theories of particle
physics act on the associated quantum theory, it was discovered
that they might not:  there could an ``anomaly'' in the quantum
theory which prevented the action from being well-defined.
This phenomenon was originally expressed in terms of Feynman diagrams,
and came to be regarded as the statement that, although the symmetry
group preserves the Lagrangian, it fails to preserve the measure
on the space of paths.

This formulation is unsatisfying to mathematicians, who know that the
measure hasn't been properly defined in mathematics and so who rightfully
wonder how an ill-defined thing can fail to be preserved by a group
action?

There is another interpretation of the anomaly, however, in terms of
the Dirac operators of the quantum theory.  The equations of motion for
the non-self-dual bosonic fields in a (Wick-rotated)
physical theory involve the Laplacian $\Delta:\mathcal{V}\to\mathcal{V}$
acting on a space of functions or differential forms on spacetimes
(sometimes bundle-valued).
However, for fermionic fields in the theory, the equations of motion involve
a Dirac operator 
$\nabla:\mathcal{V}\to\mathcal{W}$ 
which is a ``square root'' of the Laplacian
and typically does not map the classical space of spinors to itself.
Quillen \cite{quillen-det} 
introduced a  ``determinant line bundle'' associated to
the Dirac operator $\nabla$, and Bismut and Freed
\cite{MR0853982} equipped it with a connection.  The theory is 
anomaly-free if the line bundle is trivial and the connection is flat.
Parallel developments also
occured in Moscow, where Faddeev and Shatashvili \cite{Fadeev-Shatashvili} gave a formulation in terms of group cohomology.

The beauty of this approach is the connection to the Atiyah--Singer
index theorem, as mentioned in the quote from \cite{MR878768} at
the head of this section.  The original  Atiyah--Singer theorem
can be used to interpret anomalies of abelian group actions, but
the anomalies of non-abelian group actions (and gravitational anomalies)
generally require the families index theorem.
In all cases one is calculating topological obstructions to the
triviality of the determinant bundle and the connection on it.

\section{Donaldson invariants} \label{sec:donaldson}

\epigraph{The surprise produced by Donaldson's result was accentuated by the 
fact that
his methods were completely new and were borrowed from theoretical physics,
in the form of the Yang-Mills equations.}{Michael Atiyah (1986)}

As discussed in Section~\ref{sec:instantons},
the idea of Yang--Mills instantons was able to move from physics to
mathematics thanks to renewed communication between mathematicians
and physicists beginning in the mid 1970s.  
The next step, however, was truly remarkable.  In the early 1980s,
Simon Donaldson studied the Yang--Mills instanton equations on 
arbitrary compact
four-manifolds, and using them, was able to make very unexpected progress
in the study of differentiable four-manifolds \cite{MR710056}.  

For a fixed principal $G$-bundle $P$ with instanton number $k$
over a four-manifold $X$, where
$G$ is a compact Lie group with Lie algebra $\mathfrak{g}$, Donaldson 
considered the set $\mathcal{A}$ of connections on $P$
(which is an affine space acted upon by the space $\Omega^1(X,\mathfrak{g})$
of $\mathfrak{g}$-valued $1$-forms on $X$ and can be
regarded as a torsor for $\Omega^1(X, \operatorname{ad}(P))$), modulo 
the automorphism group $\mathcal{G}$ of $P$.  On the orbit space
$\mathcal{A}/\mathcal{G}$ the self-dual connections can be identified
as the kernel of the operator 
$d_A^-:\Omega^1(X,\mathfrak{g})\to \Omega^2_-(X,\mathfrak{g})$
which is the composition of covariant differentiation with projection to the
anti-self-dual part.  In his early work, Donaldson studied the
kernel of $d_A^-$ on $\mathcal{A}/\mathcal{G}$
in the case $G=SU(2)$ and $k=1$, obtaining a moduli space
$\mathcal{M}_1\subset \mathcal{A}/\mathcal{G}$ of self-dual
connections.\footnote{More precisely, the relevant space is given by
the kernel of a small perturbation of the operator $d_A^-$, which makes
the space more regular.}  He found that if the intersection form on $X$
is negative definite, then the moduli space $\mathcal{M}_1$
is a $5$-manifold away from
a finite collection of singularities corresponding to reducible connections
(i.e., connections compatible with a decomposition of the 
associated vector bundle
$E = P \otimes_{SU(2)} \mathbb{C}^2$ into a sum of two line bundles
$L\oplus L^{-1}$).
Analyzing these singularities carefully led to restrictions
on the intersection form on second cohomology,
and in particular showed that if the intersection form 
is positive  definite, then it is diagonalizable.  
This
had many remarkable consequences including the failure of the smooth version
of the h-cobordism conjecture in dimension four as well as the existence
of an ``exotic'' differentiable structure on (topological) $\mathbb{R}^4$.
In his proofs, Donaldson relied on earlier work  of Taubes \cite{MR658473}
and Uhlenbeck \cite{MR648355,MR648356}.
The significance of these results are dramatically strengthened when they are combined with Freedman's near-simultaneous solution of the topological Poincar\'e conjecture in dimension 4 \cite{Freedman}.

This work led to the award of a Fields medal to Donaldson at the 1986 ICM,
during which Donaldson's work
was presented by Atiyah \cite{MR934209}.  
As Atiyah emphasized in the quote at the head
of this section, the input from physics was one of the most remarkable
aspects of the work.  The conversation which had begun with the ADHM
solution of a problem in physics had now become a two-way conversation!

Donaldson's work on four-manifolds did not end with 
the awarding of the Fields medal.  He extended the work in a number
of directions, including a definition of
polynomial invariants on the cohomology of $X$ of arbitrary
degree, based on the moduli space $\mathcal{M}_k$ of instantons
with second Chern class $k$ \cite{MR1066174}.

\section{Topological quantum field theory} \label{sec:tqft}

\epigraph{[Witten's paper \cite{MR953828}], which introduces TQFT in the 
context of Donaldson's theory of
4-manifolds and Floer's theory of 3-manifolds, could well emerge as one of 
the most
significant works in late 20th century topology.}{Daniel Freed (1988)}

The problem which Donaldson solved was not really a problem from physics:
it was a problem from mathematics whose techniques were inspired by
physics.  There were many developments stemming from Donaldson's original
work, including a related theory in three dimensions developed by Andreas
Floer.

A good starting point for
Floer's theory is the Chern--Simons functional,
which on a three-manifold $Y$ equipped with a connection $A$ 
on a $G$-bundle over $Y$ is
the quantity
\begin{equation}
\operatorname{CS}(A) = \frac1{4\pi} \int_Y \operatorname{tr} (A \wedge F 
+ \frac23 A \wedge A \wedge A) \in \mathbb{R}/\mathbb{Z} .
\end{equation}
Equivalently, if $A_0$ is the trivial connection on the trivial 
$G$-bundle over $Y$ and $A_t=(1-t)A+tA_0$ then
\begin{equation}
\operatorname{CS}(A) = \frac1{8\pi^2} \int_{Y\times [0,1]} 
\operatorname{tr} (F \wedge  F).
\end{equation}
In its interpolation between $A$ and $A_0$ at the two ends of
$Y\times [0,1]$, the connection $A$ is an {\em instanton}\/ in the
sense the word is used in physics (see 
\cite{Coleman-uses}),
albeit an instanton of
Euclidean signature.

Floer used a small perturbation of
the function $\operatorname{CS}$ as
a Morse function on the space of connections $\mathcal{A}$.  
Although $\mathcal{A}$ is infinite dimensional and
the formal Hessian of $\operatorname{CS}$ at a critical point has
infinite sets of both positive and negative eigenvalues, 
it is of {\em Dirac type}\/ and  essentially coincides
with the operator $\star d$ acting
on $\Omega^1/d\Omega^0$, suitably extended to Lie algebra  valued forms.
What is well-defined is the {\em index}, the difference between
the positive and negative eigenvalues.  Moreover, as in Witten's
earlier interpretation of Morse theory  in terms
of quantum tunneling \cite{MR683171}, given two critical points $P$ and $Q$ the
solutions to $dA/dt =-\operatorname{grad} \operatorname{CS}$ 
on $Y\times \mathbb{R}$ with
connection $A_P$ as $t\to-\infty$ and connection $A_Q$ as
$t\to\infty$ are identified with boundary operators in a chain complex.
The resulting homology theory (which is only defined mod 8 
for $G=SU(2)$ due to the
index theorem in 4 dimensions) is
Floer homology \cite{MR956166}.

The critical points of (slightly perturbed) $\operatorname{CS}$ 
are identified with 
irreducible representations $\pi_1(Y)\to G$ and account for the
connection between Floer homology and the Casson invariant (which
counts such representations, with appropriate signs).  
Moreover, the Floer theory provides
a natural setting for Donaldson theory on four-manifolds with boundary.
including Donaldson's polynomial invariants.

Atiyah \cite{MR974342}
put Donaldson's and Floer's work together (also combining them
with some ideas about Heegard splitting to extend the theory to dimension two),
obtaining a non-relativistic quantum field theory.  Witten \cite{MR953828} went
one step further, and found the proper physical setting for the
work of Donaldson and Floer.

The relativistic quantum field theory which describes Donaldson's and
Floer's results is a topological twist of the 
usual\footnote{More precisely, the supersymmetric quantum field
theory with twice the minimal amount of supersymmetry.} supersymmetric
quantum field theory associated to the $SU(2)$ gauge group,
a new notion which Witten introduced in order to provide the physical
setting.  Starting from a supersymmetric theory with certain
supercharges, the action of those supercharges on the physical fields
was modified in a way whch made them independent of the choice
of metric on spacetime.  The correlation functions in the corresponding
field theory turned out to precisely be the Donaldson polynomial invariants!

Witten used path integrals to motivate his construction, and
the theory itself clearly belongs to physics, not mathematics.
But as Dan Freed remarks in the review of \cite{MR953828} quoted at the
head of this section \cite{tqft-review}, the implications for topology
itself were profound.  In fact, the intrinsic study of topological
field theories (independent of the precise details of Donaldson theory)
has become an important aspect of twenty-first century mathematics.

\section{Seiberg--Witten theory} \label{sec:seiberg-witten}

\epigraph{In the last three months of 1994 a remarkable thing happened:
this research area was turned on its head by the introduction of a new
kind of differential-geometric equation by Seiberg and Witten:  in the space
of a few weeks long-standing problems were solved, new and unexpected results
were found, along with simpler new proofs of existing ones, and new vistas
for research opened up.}{Simon Donaldson (1995)}

In 1994, Seiberg and Witten made some of the first progress in understanding
quantum field theories from a non-perturbative perspective (in particular, studying
properties which are not necessarily tied to the 
path-integral formalism \cite{Seiberg:1994rs,MR1303306,Seiberg:1994aj}).  
The particular
theory which they studied first -- the $\mathcal{N}=2$ supersymmetric
gauge theory in four dimensions -- was the same one which Witten had
earlier shown could be topologically twisted  to yield the Donaldson--Witten
theory.  Now, Seiberg and Witten were able to find a new description
of the infrared behavior of the $SU(2)$ gauge theory which took the form of
a $U(1)$ gauge theory coupled to a magnetic monopole.

This work in physics had an immediate consequence in mathematics 
(after twisting) in the
form of new topological invariants analogous to the Donaldson invariants
 \cite{Witten:1994cga}.
These Seiberg--Witten invariants, as they came to be called, were 
substantially easier to compute than the Donaldson invariants and
progress was quickly  made on many difficult conjectures which had
been left open by the original Donaldson theory.  Donaldson himself
expressed great astonishment at the speed of progress 
\cite{MR1339810}, as quoted at the head of this section.

As John Morgan said during a lecture at the 1995 Cornell Topology Festival
\cite{Morgan:unpublished}: ``The physicists keep coming up with amazing 
equations for us to solve.  Once we know the equations, we can get lots of 
mathematics out of them, but why can't we find the equations ourselves?''

The exchange of ideas between physics and mathematics had now proceeded
through at least four stages:
from progress in solving the instanton 
equation in physics, to Donaldson's application of those ideas to
the understanding to four-manifolds, to Witten's construction of topological
quantum field theory, to the Seiberg--Witten study of the infrared properties
of that theory and the spectacular mathematics which resulted!

\section{Conclusions} \label{sec:conclusions}

\epigraph{I am sure that [the interaction of math and physics] 
is going to continue and I
believe the reason it will continue is that quantum
field theory and string theory \dots\ have
rich mathematical secrets.}{Edward Witten (2014)}

The story I have presented about interactions between mathematics
and physics did not end in 1994, but has continued to
develop fruitfully in many directions.  For lack of time, I did not mention 
string
theory at all in this lecture, but many of the important interactions since
the mid 1980s have
involved string theory as well as quantum field theory.  
As Edward
Witten predicted in a 2014 interview conducted by 
Hirosi Ooguri \cite{MR3328919}
and quoted at the head of this section,
this interaction is likely to continue for a long time to come!
Witten went on to say:  ``When some of
these secrets come to the surface, they often come
as surprises to physicists because we do not really
understand string theory properly as physics -- we
do not understand the core ideas behind it. At an
even more basic level, the mathematicians are still
not able to fully come to grips with quantum field
theory and therefore things coming from it are
surprises. So for both of those reasons, I think that
the physics and math ideas generated are going to
be surprising for a long time.''


\begin{thebibliography}{AHDM78}

\bibitem[AG86]{MR878768}
Luis Alvarez-Gaum{\'e}, \emph{An introduction to anomalies}, Fundamental
  problems of gauge field theory ({E}rice, 1985), NATO Adv. Sci. Inst. Ser. B
  Phys., vol. 141, Plenum, New York, 1986, pp.~93--206. \MR{878768 (88c:81077)}

\bibitem[AHDM78]{MR598562}
M.~F. Atiyah, N.~J. Hitchin, V.~G. Drinfel{\cprime}d, and Yu.~I. Manin,
  \emph{Construction of instantons}, Phys. Lett. A \textbf{65} (1978), no.~3,
  185--187. \MR{598562}

\bibitem[AHS77]{MR0458424}
M.~F. Atiyah, N.~J. Hitchin, and I.~M. Singer, \emph{Deformations of
  instantons}, Proc. Nat. Acad. Sci. U.S.A. \textbf{74} (1977), no.~7,
  2662--2663. \MR{0458424}

\bibitem[Ati79]{MR554924}
M.~F. Atiyah, \emph{Geometry of {Y}ang-{M}ills fields}, Scuola Normale
  Superiore Pisa, Pisa, 1979. \MR{554924}

\bibitem[Ati87]{MR934209}
Michael Atiyah, \emph{On the work of {S}imon {D}onaldson}, Proceedings of the
  {I}nternational {C}ongress of {M}athematicians, {V}ol. 1, 2 ({B}erkeley,
  {C}alif., 1986), Amer. Math. Soc., Providence, RI, 1987, pp.~3--6.
  \MR{934209}

\bibitem[Ati88]{MR974342}
\bysame, \emph{New invariants of {$3$}- and {$4$}-dimensional manifolds}, The
  mathematical heritage of {H}ermann {W}eyl ({D}urham, {NC}, 1987), Proc.
  Sympos. Pure Math., vol.~48, Amer. Math. Soc., Providence, RI, 1988,
  pp.~285--299. \MR{MR974342 (89m:57034)}

\bibitem[BF86]{MR0853982}
Jean-Michel Bismut and Daniel~S. Freed, \emph{The analysis of elliptic
  families. {I}. {M}etrics and connections on determinant bundles}, Comm. Math.
  Phys. \textbf{106} (1986), no.~1, 159--176. \MR{853982 (88h:58110a)}

\bibitem[BH10]{MR2651086}
John Baez and John Huerta, \emph{The algebra of grand unified theories}, Bull.
  Amer. Math. Soc. (N.S.) \textbf{47} (2010), no.~3, 483--552. \MR{2651086
  (2012b:81271)}

\bibitem[BPST75]{MR0434183}
A.~A. Belavin, A.~M. Polyakov, A.~S. Schwartz, and Yu.~S. Tyupkin,
  \emph{Pseudoparticle solutions of the {Y}ang-{M}ills equations}, Phys. Lett.
  B \textbf{59} (1975), no.~1, 85--87. \MR{0434183}

\bibitem[Col85]{Coleman-uses}
Sidney~R. Coleman, \emph{The uses of instantons}, Aspects of Symmetry,
  Cambridge University Press, 1985, pp.~265--350.

\bibitem[Dir31]{dirac-quantised}
P.~A.~M. Dirac, \emph{Quantised singularities in the electromagnetic field},
  Proceedings of the Royal Society of London A: Mathematical, Physical and
  Engineering Sciences \textbf{133} (1931), no.~821, 60--72.

\bibitem[Don83]{MR710056}
S.~K. Donaldson, \emph{An application of gauge theory to four-dimensional
  topology}, J. Differential Geom. \textbf{18} (1983), no.~2, 279--315.
  \MR{MR710056 (85c:57015)}

\bibitem[Don90]{MR1066174}
\bysame, \emph{Polynomial invariants for smooth four-manifolds}, Topology
  \textbf{29} (1990), no.~3, 257--315. \MR{MR1066174 (92a:57035)}

\bibitem[Don96]{MR1339810}
\bysame, \emph{The {S}eiberg-{W}itten equations and {$4$}-manifold topology},
  Bull. Amer. Math. Soc. (N.S.) \textbf{33} (1996), no.~1, 45--70. \MR{1339810}

\bibitem[Dys72]{MR0522147}
Freeman~J. Dyson, \emph{Missed opportunities}, Bull. Amer. Math. Soc.
  \textbf{78} (1972), 635--652. \MR{0522147 (58 \#25442)}

\bibitem[FS84]{Fadeev-Shatashvili}
L. D. Faddeev and S. L. Shatashvili, \emph{Algebraic and Hamiltonian
methods in the theory of non-Abelian anomalies}, Theoret. and Math. Phys.
\textbf{60} (1984), no.~2, 770-778, DOI 10.1007/BF01018976.


\bibitem[Flo88]{MR956166}
Andreas Floer, \emph{An instanton-invariant for {$3$}-manifolds}, Comm. Math.
  Phys. \textbf{118} (1988), no.~2, 215--240. \MR{956166}

\bibitem[Fre89]{tqft-review}
Daniel~S. Freed, \emph{Review of \cite{MR953828}}, MR953828, 1989.

\bibitem[Fre82]{Freedman}
Michael H. Freedman, The topology of four-dimensional manifolds,
J. Diff Geom. 17 (1982) 357--453.

\bibitem[Lic]{MR0088015}
Andr\'{e} Lichnerowicz, \emph{Th\'{e}orie globale des connexions et des groupes
  d'holonomie}, Edizioni Cremonese, Roma, 1957. \MR{0088015}

\bibitem[Mor]{Morgan:unpublished}
John Morgan, \emph{Applications of {S}eiberg--{W}itten invariants to the
  topology of 4-manifolds}, Lecture at the 1995 Cornell Topology Festival,
  unpublished.

\bibitem[Oog15]{MR3328919}
Hirosi Ooguri, \emph{Interview with {E}dward {W}itten}, Notices Amer. Math.
  Soc. \textbf{62} (2015), no.~5, 491--506. \MR{3328919}

\bibitem[Pol75]{Polyakov:1975rs}
A.~M. Polyakov, \emph{Compact gauge fields and the infrared catastrophe}, Phys.
  Lett. \textbf{B59} (1975), 82--84.

\bibitem[PS95]{MR1402248}
Michael~E. Peskin and Daniel~V. Schroeder, \emph{An introduction to quantum
  field theory}, Addison-Wesley Publishing Company Advanced Book Program,
  Reading, MA, 1995, Edited and with a foreword by David Pines. \MR{1402248
  (97j:81001)}

\bibitem[Qui85]{quillen-det}
D.~Quillen, \emph{Determinants of {C}auchy-{R}iemann operators on {R}iemann
  surfaces}, Funktsional. Anal. i Prilozhen. \textbf{19} (1985), no.~1, 37--41,
  96, English translation: Functional Anal. Appl. {\bf 19} (1985), no. 1,
  31--34. \MR{783704}

\bibitem[Ste51]{MR0039258}
Norman Steenrod, \emph{The {T}opology of {F}ibre {B}undles}, Princeton
  Mathematical Series, vol. 14, Princeton University Press, Princeton, N. J.,
  1951. \MR{0039258}

\bibitem[SW94a]{Seiberg:1994rs}
N.~Seiberg and E.~Witten, \emph{Electric-magnetic duality, monopole
  condensation, and confinement in {$N=2$} supersymmetric {Y}ang-{M}ills
  theory}, Nuclear Phys. B \textbf{426} (1994), no.~1, 19--52. \MR{1293681}

\bibitem[SW94b]{MR1303306}
\bysame, \emph{Erratum: ``{E}lectric-magnetic duality, monopole condensation,
  and confinement in {$N=2$} supersymmetric {Y}ang-{M}ills theory''}, Nuclear
  Phys. B \textbf{430} (1994), no.~2, 485--486. \MR{1303306}

\bibitem[SW94c]{Seiberg:1994aj}
\bysame, \emph{Monopoles, duality and chiral symmetry breaking in {$N=2$}
  supersymmetric {QCD}}, Nuclear Phys. B \textbf{431} (1994), no.~3, 484--550.
  \MR{1306869}

\bibitem[TAASSSS]{MR3748689}
L.~A. Takhtadzhyan, A.~Yu. Alekseev, I.~Ya. Aref\cprime~eva, M.~A.
  Semenov-Tyan-Shanski\u{i}, E.~K. Sklyanin, F.~A. Smirnov, and S.~L.
  Shatashvili, \emph{L. {D}. {F}addeev's scientific legacy. {A} survey of his
  work}, Uspekhi Mat. Nauk \textbf{72} (2017), no.~6(438), 3--112. \MR{3748689}

\bibitem[Tau82]{MR658473}
Clifford~Henry Taubes, \emph{Self-dual {Y}ang-{M}ills connections on
  non-self-dual {$4$}-manifolds}, J. Differential Geom. \textbf{17} (1982),
  no.~1, 139--170. \MR{658473}

\bibitem[Uhl82a]{MR648355}
Karen~K. Uhlenbeck,
\emph{Removable singularities in {Y}ang-{M}ills fields}, Comm. Math.
  Phys. \textbf{83} (1982), no.~1, 11--29. \MR{648355}

\bibitem[Uhl82b]{MR648356}
\bysame,
\emph{Connections with {$L^{p}$}\ bounds on curvature},
  Comm. Math. Phys. \textbf{83} (1982), no.~1, 31--42. \MR{648356}

\bibitem[Wig60]{wigner-unreasonable}
Eugene~P. Wigner, \emph{The unreasonable effectiveness of mathematics in the
  natural sciences}, Comm. Pure Appl. Math. \textbf{13} (1960), 1--14.

\bibitem[Wit82]{MR683171}
Edward Witten, \emph{Supersymmetry and {M}orse theory}, J. Differential Geom.
  \textbf{17} (1982), no.~4, 661--692 (1983). \MR{MR683171 (84b:58111)}

\bibitem[Wit88]{MR953828}
\bysame, \emph{Topological quantum field theory}, Comm. Math. Phys.
  \textbf{117} (1988), no.~3, 353--386. \MR{953828}

\bibitem[Wit94]{Witten:1994cga}
\bysame, \emph{Monopoles and four-manifolds}, Math. Res. Lett. \textbf{1}
  (1994), no.~6, 769--796. \MR{1306021}

\bibitem[WY75]{MR0426712}
Tai~Tsun Wu and Chen~Ning Yang, \emph{Concept of nonintegrable phase factors
  and global formulation of gauge fields}, Phys. Rev. D (3) \textbf{12} (1975),
  no.~12, 3845--3857. \MR{0426712}

\bibitem[YM54]{MR0065437}
C.~N. Yang and R.~L. Mills, \emph{Conservation of isotopic spin and isotopic
  gauge invariance}, Physical Rev. (2) \textbf{96} (1954), 191--195.
  \MR{0065437}

\bibitem[Zim08]{simons-yang}
Bill Zimmerman, \emph{Interview with {C}hen {N}ing {Y}ang and {J}ames {H}arris
  {S}imons}, Stony Brook Masters Series, March 28, 2008, Available at
  https://www.youtube.com/watch?v=zVWlapujbfo.

\end{thebibliography}

\providecommand{\bysame}{\leavevmode\hbox to3em{\hrulefill}\thinspace}
\providecommand{\MR}{\relax\ifhmode\unskip\space\fi MR }
\providecommand{\MRhref}[2]{%
  \href{http://www.ams.org/mathscinet-getitem?mr=#1}{#2}
}
\providecommand{\href}[2]{#2}

\end{document}